\documentclass[a4paper,12pt,fleqn]{article}

\usepackage{amssymb}
\usepackage{amsmath}
\usepackage{graphicx}
\usepackage{color}
\usepackage{eurosym}
\usepackage{txfonts}
\usepackage{courier}
\usepackage{pifont}

\usepackage[colorlinks=false,urlcolor=blue,ps2pdf=true,
pdfpagemode=UseNone, pdfauthor={Robin Whitty},
pdfsubject={Mathematics of OR}, pdfcreator={Robin Whitty},
pdfproducer={Robin Whitty}]{hyperref}

\newcommand{\ds}{\displaystyle}
\newcommand{\bm}[1]{\mbox{\boldmath $#1$}}

\newcommand{\QED}{\unskip\nobreak\hfill\rlap{\llap{\large $\Box$}}} 

\newtheorem{Thm}{Theorem}

\newtheorem{Exmpl}[Thm]{Example}

\newtheorem{Lm}[Thm]{Lemma}
\newtheorem{Cr}[Thm]{Corollary}

\newtheorem{Rem}[Thm]{Remarks}

\author{Fatima Affif Chaouche\thanks{University of Sciences and Technology Houari Boumediene, Algiers},
Carrie Rutherford\thanks{London South Bank University}, Robin
Whitty\thanks{Queen Mary University of London}}
\date{\today}
\title{Pancyclicity when each cycle must pass exactly $k$ Hamilton cycle chords}

\begin{document}

\maketitle

\hspace{.25in}
\parbox[c]{6in}{\small {\bf Abstract} It is known that $\Theta(\log n)$ chords must be added to an $n$-cycle to produce a pancyclic graph; for vertex pancyclicity, where every vertex belongs to a cycle of every length, $\Theta(n)$ chords are required. A possibly `intermediate' variation is the following: given $k$, $1\leq k\leq n$, how many chords must be added to ensure that there exist cycles of every possible length each of which passes exactly $k$ chords? For fixed $k$, we establish a lower bound of $\Omega\big(n^{1/k}\big)$ on the growth rate.
}


A simple graph $G$ on $n$ vertices is {\em pancyclic} if it has cycles of every length $l$, $3\leq l\leq n$. The study of these graphs was initiated by Bondy's observation \cite{Bondy, Bondy2} that, for non-bipartite graphs, sufficient conditions for hamiltonicity can also be sufficient for pancyclicity.
In general, we may distinguish, in a pancyclic graph $G$, a Hamilton cycle $C$; then the remaining edges of $G$ form chords  of $C$. We can then ask, given $k\leq l\leq n$ if, relative to $C$, a cycle of length $l$ exists which uses exactly $k$ chords. This suggests a $k$-chord analog of pancyclicity: do all possible cycle lengths occur when cycles must use exactly $k$-chords of a suitably chosen Hamilton cycle?

We accordingly define a function $c(n,k)$, $n\geq 6,$ $k\geq 1$, to be the smallest number of chords which must be added to an $n$-cycle in order that cycles of all possible lengths may be found, each passing exactly $k$ chords.  No Hamilton cycle can use exactly one chord of another Hamilton cycle, so that when $k=1$ cycle lengths must lie between $k$ and $n-1$. The function is undefined for $k>n$. We define the function for $n\geq 6$ because $n=4,5$ are too restrictive to be of interest to us.

Our aim in this paper is to investigate the growth of the function $c(n,k)$ as $n$ increases, for fixed $k$. 

\begin{Exmpl}
\label{exmpl:pancycles}
 Label the vertices around the cycle $C_6$, in order, as $v_1,\ldots,v_6$. Add chords $v_1v_3$ and $v_1v_4$; the result is a pancyclic graph. It also has cycles of all lengths $\leq 5$ each passing exactly one of the chords. If $v_2v_6$ is added then cycles exist of all lengths $\geq 3$, each passing two chords. If two further chords, $v_2v_4$ and $v_4v_6$, are added then cycles exist of all lengths $\geq 3$, each passing three chords. For 4-chord cycles we require six chords to be added, i.e. $C(6,4)=6$. Six suitably chosen chords are also sufficient for $5-chord$ and $6-chord$ cycles: $C(6,5)=C(6,6)=6$.
\end{Exmpl}

\begin{Lm}
 \label{lm:basic} \begin{enumerate}
  \renewcommand{\labelenumi}{(\arabic{enumi})}
\item $\ds c(n,1)=\left\lfloor\frac{n-3}{2}\right\rfloor$.
\item $c(n,k)\geq k,$ with equality if and only if $k=n$.
\item $c(n,n-1)=n$.
\end{enumerate}
\end{Lm}

{\bf Proof.} (1) follows from the observation that a chord in $C_n$ forming a 1-chord cycle of length $k$ automatically forms a 1-chord cycle of length $n+2-k$.

(2) is immediate from the definition of $c(n,k)$.

(3) Let $G$ consist of an $(n-1)$-cycle, together with an $(n-1)$-chord cycle on the same vertices. Choose vertex $v$: let the chords at $v$ be $xv$ and $yv$ and its adjacent cycle edges be $uv$ and $vw$, with $u,v,w,x,y$ appearing in clockwise order around the cycle. Replace $v$ and its incident edges with two vertices $v_u$ and $v_w$, with edges $v_uv_w$, $uv_u$, $v_ww$, $xv_w$ and $yv_u$. The $(n-1)$-chord cycle in $G$ becomes an $(n-1)$-chord $n$-cycle. Add an $n$-th chord $xv_u$ to give an $(n-1)$-chord $(n-1)$-cycle.\QED

Table 1 supplies some small values/bounds for $c(n,k)$. The lower bounds are supplied by Corollary~\ref{cr:pancyclic1} (see below); except for those values covered by Lemma~\ref{lm:basic}, exact values and upper bounds were found by computer search.

$\begin{array}{rr|cccccccccccccccccccccc}
      &       && \bm{k} &&\\
      &       && \bm{1} && \bm{2} && \bm{3} && \bm{4} && \bm{5} && \bm{6} && \bm{7} && \bm{8} && \bm{9} && \bm{10}  && \bm{11}
      \\\hline
 &       &\\[-.05in]
\bm{n}& \bm{6} && 2 && 3 && 5      && 6       && 6      && 6\\[.1in]
      & \bm{7} && 2 && 3 && 5      && 6       && 6      && 7 && 7\\[.1in]
      & \bm{8} && 3 && 4 && 5      && 6       && 6      && 7 && 8 && 8\\[.1in]
      & \bm{9} && 3 && 4 && 5      && 6       && 7      && 8  && 8  && 9 && 9\\[.1in]
      &\bm{10} && 4 && 4 && 5      && 6       && \geq 6 && \geq 7  && \geq 8  && \geq 9 && 10 && 10\\[.1in]
      &\bm{11} && 4 && 4 && \geq 5 && \geq 6  && \geq 7 && \geq 7  && \geq 8  && \geq 9  && \geq 10  && 11 && 11\\[.1in]
      &\bm{12} && 5 && 4 && \geq 5 && \geq 6  && \geq 7 && \geq 7  && \geq 8  && \geq 9  && \geq 10 &&  \geq 11  && 12  \\[.1in]
      &\bm{13} && 5 && 4 && \geq 5 && \geq 6  && \geq 7 && \geq 8  && \geq 8  && \geq 9  && \geq 10 &&   \geq 11 && \geq 12     \\[.1in]
\end{array}$

\centerline{Table 1. Values of $c(n,k)$ for $6\leq n\leq 13$ and $1\leq k\leq 11$.}

Our aim is to compare $c(n,k)$ with the number of chords required for pancyclicity and for {\em vertex pancyclicity}, in which each vertex must lie on a cycle of every length.

The following lower bound is stated without proof in \cite{Bondy}:

\begin{Thm}
\label{thm:pancyclic} In a pancyclic graph $G$ on $n$ vertices the number of edges is not less than $n-1+log_2(n-1)$.\QED
\end{Thm}

For the sake of completeness we observe that theorem~\ref{thm:pancyclic} follows immediately from the following lemma:

\begin{Lm}
 \label{lm:pancyclic}Suppose $p$ chords are added to $C_n$, $n\geq 3$. Then the number $N(n,p)$ of cycles in the resulting graph satisfies
$${p+2 \choose 2}\leq N(n,p)\leq 2^{p+1}-1.$$
\end{Lm}

{\bf Proof.} Embed $C_n$ convexly in the plane. Suppose the chords added to $C_n$ are, in order of inclusion, $e_1, e_2,\ldots, e_p$. Say that $e_i$ intersects $e_j$ if these edges cross each other when added to the embedding of $C_n$. Let $n_i$ be the number of new cycles obtained with $e_i$ is added. Then $n_i$ satisfies:
\begin{enumerate}
\item $n_i\geq i+1$, the minimum occurring if and only if the $e_j$ are pairwise non-intersecting for  $j\leq i$;
\item $n_i\leq 2^i$, the maximum occurring if and only if $e_i$ intersects with $e_j$ for all $j<i$, giving
$n_i=\sum_{j=0}^{i}{i\choose j}$.
\end{enumerate}
Now $\ds1+\sum_{i=1}^p (i+1)\leq 1+\sum_{i=1}^p n_i \leq 1+\sum_{i=1}^p 2^i$ and the result follows.
 \QED

The exact value of the minimum number of edges in an $n$-vertex pancyclic graph has been calculated for small $n$ by George et al \cite{George} and Griffin \cite{Griffin}. For $3\leq n\leq 14$ the lower bound in theorem~\ref{thm:pancyclic} is exact; however, it can be seen that, for $n=15, 16$,  we must add four chords to $C_n$ to achieve pancyclicity while the argument in the proof of lemma~\ref{lm:pancyclic} can only account for three.

As regards an upper bound on the number of chords required for pancyclicity, \cite{Bondy} again asserts $O(\log n)$, again without a proof. A $\log n$ construction has been given by Sridharan \cite{Sridharan}. Together with theorem~\ref{thm:pancyclic} this gives an `exact' growth rate for pancyclicity: it is achieved by adding $\Theta(\log n)$ chords to~$C_n$.

In contrast, {\em vertex pancyclicity}, in which every vertex lies in a cycle of every length has been shown by Broersma \cite{Broersma} to require $\Theta(n)$ edges to be added to $C_n$. Our question is: where between $\log n$  and $n$ does $c(n,k)$ lie? For fixed $k$, we find a lower bound strictly between the two: $\Omega(n^{1/k})$.

Let us for the moment restrict to $k\geq 3$. Suppose we add $p$ chords to $C_n$, $\ds 3\leq k\leq p\leq {n\choose 2}-n$. Suppose that these $p$ added chords include a $k$-cycle. We will use $K(k,p)$, defined for $k\geq 3$, to denote the maximum number of $k$-chord cycles that can be created in the resulting graph. Then $1\leq K(k,p)$ by definition and $K(k,p)\leq 2^{p+1}-1$ by lemma~\ref{lm:pancyclic}. By lowering this upper bound we can increase the lower bound on C(n,k).

\begin{Thm}
\label{thm:kchords} $\ds K(k,p)\leq {p\choose k}+k{p-k\choose k-1}+{p-k\choose k}.$
\end{Thm}

We will use the following Lemma to prove theorem~\ref{thm:kchords}:
\begin{Lm}
\label{lm:addchord}
Suppose that a set $X$ of  chords is added to $C_n$. In the resulting graph the maximum number of 
cycles passing all edges in $X$ is
$$\left\{\begin{array}{ccl}
1 &&\mbox{if $X$ contains adjacent chords}\\
2&&\mbox{if no two chords of $X$ are adjacent}
\end{array}\right.$$
\end{Lm}
{\bf Proof.} Let $G$ be the graph resulting from adding the chords of $X$ to $C_n$. We may assume without loss of generality that $G$ has no vertices of degree 2, since such vertices may be contracted out. For a given cycle in $G$ passing all chords of $X$, let $H$ denote the intersection of this cycle with the $C_n$. Then $H$ consists of isolated vertices and disjoint edges, and $H$ is completely determined once any of these vertices or edges is fixed. If two chords are adjacent this fixes an isolated vertex of $H$; if no two chords are adjacent then there is a maximum of two ways in which a single edge of $H$ may be fixed.
\QED

{\bf Proof of theorem~\ref{thm:kchords}.} By definition of $K(k,p)$ we must use a set, say $S$, of $k$ chords to create a $k$-cycle. We add new chords to $S$, one by one. On adding the $r$-th additional chord, $1\leq r\leq p-k$, we ask how many $k$-chord cycles use this chord. For any such a cycle the previous $r-1$ chords will be split between $S$ and non-$S$ chords: with $i$ chords from $S$ being used, $0\leq i\leq k-1$, this can happen in
$${k\choose i}{r-1\choose k-i-1}$$
ways. Since $i>1$ forces two adjacent chords in $S$ to be used, summing over $i$, according to lemma~\ref{lm:addchord}, and then over~$r$ gives
$$K(k,p)\leq 1+\sum_{r=1}^{p-k}\left(2\sum_{i=0}^1{k\choose i}{r-1\choose k-i-1}
+\sum_{i=2}^{k-1}{k\choose i}{r-1\choose k-i-1}\right).$$
This simplifies (e.g. using symbolic algebra software such as Maple) to give the result.
\QED

\begin{Cr}
\label{cr:pancyclic1}
For given positive integers $k$ and $n$, with $3\leq k\leq n$ and $n\geq 6$, the value of $c(n,k)$ is not less than the largest root of the following polynomial in $p$:
$$\Pi(p;n,k)={p\choose k}+k{p-k\choose k-1}+{p-k\choose k}-n+k-1.$$\QED
\end{Cr}

We finally extend our analysis to include the cases $k=1,2$:

\begin{Cr}
\label{cr:pancyclic2}
Let $n\geq 6$ be a positive integer. Then for $k\geq 1$ fixed, $c(n,k)$ is of order $\Omega\big(n^{1/k}\big)$.
\end{Cr}

{\bf Proof.}
For $k=1$ the required linear bound was provided in lemma~\ref{lm:basic}.

For $k=2$ an analysis similar to that used in the proof of theorem~\ref{thm:kchords} shows that the number of $2$-chord cycles which may be created by adding $p$ chords to $C_n$ is at most $p^2-p-1$. So to have $2$-chord cycles of all lengths from 3 to $n$ we require $p^2-p-1\geq n-2$. In this case we can solve explicitly to get the bound $\ds p\geq\frac12\left(1+\sqrt{4n-3}\,\right)$.

Now suppose $k\geq 3$. In order to have all $k$-chord cycles of all lengths between $k$ and $n$ we must have
$$n-k+1\leq {p\choose k}+k{p-k\choose k-1}+{p-k\choose k}\leq f(k)p^k,$$
for some function $f(k)$. Therefore $p^k\geq (n-k+1)/f(k)$ so, for $k$ fixed, $p=\Omega\big(n^{1/k}\big)$.
\QED

\begin{Rem}
\begin{enumerate}
\item
We are suggesting that the value of $c(n,k)$ may be `intermediate' between pancyclicity and vertex pancyclicity in the sense that the number of chords it requires to be added to $C_n$ may lie between $\log n$ and $n$. Thus far we have only a lower bound in support of our suggestion. Moreover, a comparison of the growth orders, $\Omega(\log n)$ as opposed to $\Omega\big(n^{1/k}\big)$, suggests that this is very much a `for large $n$' type result. The equation $\ln n=n^{1/k}$ has two positive real solutions for $k\geq 3$, given in terms of the two real branches of the Lambert $W$ function \cite{Corless}. In particular $\ln n$ exceeds $n^{1/k}$ for $n>e^{-kW_{-1}(-1/k)}$, and this bound grows very fast with $k$: at least  two orders of magnitude per unit increase! To give a specific example, $k=10$, the $\log$ bound exceeds the 10-th root bound until the number of vertices exceeds about $3.4\times 10^{15}$. Until then, so far as our analysis goes, we might expect `most' pancyclic graphs to be 10-chord pancyclic. However we suggest that, in the long term, a guarantee of this implication, analogous to hamiltonicity guaranteeing pancyclicity, will not be found.

\item We would like to know if $c(n,k)$ is monotonically increasing in $n$. However, it is still open even whether pancyclicity is monotonic in the number of chords requiring to be added to $C_n$ (the question is investigated in \cite{Griffin}). We believe that $c(n,k)$ it is not increasing in~$k$ and $c(n,1)>c(n,2)$ for $n=12,13$ confirms this in a limited sense. Our $n^{1/k}$ lower bound instead suggests the possibility that $c(n,k)$ is convex for fixed $n$, as a function of~$k$.

\item  We observe that, unlike pancyclicity, the property of having cycles of all lengths each passing $k$ chords is not an invariant of a graph: it depends on the initial choice of a Hamilton cycle. For example, in  figure 1, there are cycles of all lengths $\leq 9$ each passing exactly one of the $c(10,1)=4$ chords of the outer cycle but there is no 4-cycle passing exactly one chord of the bold-edge Hamilton cycle.
 \end{enumerate}
 \end{Rem}

\begin{figure}[htb]
\centerline{\includegraphics[scale=0.30]{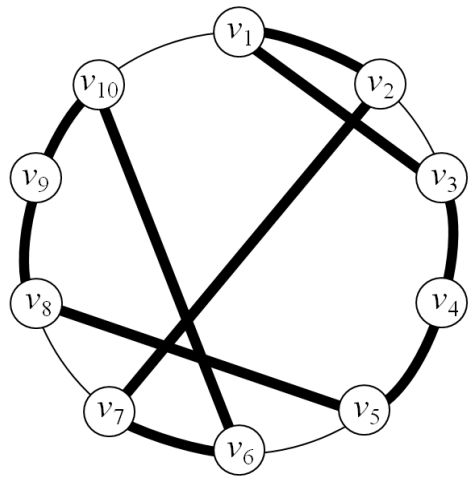}}

\caption{No 4-cycle uses exactly 1 chord of the bold-edge Hamilton cycle.}
\label{fig:cycles}
\end{figure}

\end{document}